\documentclass[11pt]{article}

\usepackage{mathrsfs}
\usepackage{latexsym}
\usepackage{amsmath,amsthm}
\usepackage{amsfonts}
\usepackage{url}

\usepackage{amssymb}

\usepackage{mytheo}
\newcommand{\atwo}{\alpha}

\renewcommand{\k}{^{(k)}}

\usepackage{bbm}
\newcommand{\pr}[1]{\P\!\paren{#1}}

\newcommand{\E}{\mathbb{E}}
\renewcommand{\P}{\mathbb{P}}

\usepackage{geometry}
\makeatletter

\newcommand{\ben}{\begin{enumerate}}
\newcommand{\een}{\end{enumerate}}
\newcommand{\bit}{\begin{itemize}}
\newcommand{\eit}{\end{itemize}}

\newcommand{\nrm}[1]{\left\Vert #1 \right\Vert}

\newcommand{\beq}{\begin{eqnarray*}}
\newcommand{\eeq}{\end{eqnarray*}}
\newcommand{\beqn}{\begin{eqnarray}}
\newcommand{\eeqn}{\end{eqnarray}}
\newcommand{\paren}[1]{\left( #1 \right)}

\newcommand{\tlprn}[1]{\left\{ #1 \right\}}
\newcommand{\set}[1]{\tlprn{#1}}
\newcommand{\abs}[1]{\left| #1 \right|}

\newcommand{\hide}[1]{}
\newcommand{\oo}[1]{\frac{1}{#1}}

\newcommand{\wed}[1]{\overset{{\mbox{\scriptsize $#1$}}}{\to}} %

\title{
The state complexity of random DFAs
}
\author{
Daniel Berend
and
Aryeh Kontorovich
}

\begin{document}
\maketitle
\begin{abstract}
The state complexity of a Deterministic Finite-state automaton (DFA)
is the number of states in its minimal equivalent DFA.
We study the state complexity of random $n$-state DFAs 
over a $k$-symbol alphabet, drawn uniformly
from the set $[n]^{[n]\times[k]}\times2^{[n]}$ of all such automata.
We show that, with high probability, the latter is 
$\alpha_k n + O(\sqrt n\log n)$
for a certain
explicit constant $\alpha_k$.
\end{abstract}

\section{Introduction}
A randomly generated deterministic 
finite
automaton (DFA) 
on $n$ states
and $k$ symbols
is drawn as follows: 
for each state and each 
of the $k$
symbols
in the alphabet,
the transition arrow's destination is chosen uniformly at random among the $n$ states; the $n
k
$
random choices are independent\footnote{By symmetry, we may always take
the state $q=1$ to be the starting state.}. 
Then each state is chosen to be accepting (or not) 
independently with probability $1/2$.
This natural model for a ``typical'' DFA 
goes back to \cite{Trakhtenbrot73} and
was considered in \cite{DBLP:conf/alt/AngluinEKR10,kont-nadler08} 
in the context of learning theory. In particular, 
in \cite{DBLP:conf/alt/AngluinEKR10} it is shown (perhaps surprisingly) 
that random DFAs possess sufficient complexity to 
embed nontrivial parity problems.

Let us define the {\em state complexity}
of a DFA $M$ as the number of states in the canonical (minimal) DFA 
equivalent to $M$, and denote it by $\nrm{M}$.
In this paper, we study the 
state complexity
of random DFAs in the model defined above.
\paragraph{Related work}
We are not aware of previous literature dealing with the specific problem we have posed.
The somewhat related problem of enumerating finite automata according to various criteria has
been extensively studied; see
\cite{MR1990452} and the references therein.
Some recent results include enumeration of minimal automata
\cite{DBLP:conf/stacs/BassinoDS12}, generation of random complete DFAs
\cite{MR2114149}, and enumeration and generation of accessible DFAs \cite{MR2347395}.
In a different line of enquiry, Pittel investigated the distributions induced by
transitive closures \cite{MR690140} and rumor spreading \cite{MR873245}.

\section{Background and notation}
We use standard automata-theoretic notation throughout; the reader is
referred to \cite{DBLP:books/daglib/0011126,ShengYu}
for background.
We put
$[n]=\set{0,\ldots,n-1}$.
Thus, $[k]$ is a $k$-ary alphabet and
$[k]^*$ 
is the 
set of all finite 
words (strings) over 
this alphabet.
The notation
$\abs{
\cdot
}$
is used for both word length and set cardinality.
Standard order-of-magnitude notation $o(\cdot)$ and $O(\cdot)$ is used,
as well as their ``with high probability'' variants $o_P(\cdot)$ 
and $O_P(\cdot)$. 
The $\tilde O(\cdot)$ notation ignores polylog factors.

An $n$-state $k$-ary Deterministic Finite-state Automaton 
is 
a
tuple
$M=(Q,q_0,A,\delta)$ where
\bit
\item $Q=[n]$
is the set of states
\item $q_0=1$ is the starting state;
\item $A\subseteq [n]$ is the set of accepting states;
\item $\delta:[n]\times[k]\to [n]$ is the transition function.
\eit
The 
transition function $\delta$ may be extended to
$[n]\times[k]^*$ via the recursion
\beqn
\label{eq:drecurs}
\delta(q,u_1u_2\cdots u_{n})
= \delta( \delta(q,u_1),u_2\cdots u_{n})
. 
\eeqn

If the accepting states are unspecified, 
the transition function $\delta$ induces
a directed multigraph on $n$ nodes with regular outdegree $k$,
called 
a $k$-ary {\em semiautomaton}. 

We recall the standard equivalence relation over the states of a DFA:
a word $x\in[k]^*$
{\em distinguishes} between
the states $p,q\in[n]$
if exactly one of 
the states $\delta(p,x)$, $\delta(q,x)$ is accepting.
If no $x\in[k]^*$ distinguishes between $p$ and $q$,
these states are {\em equivalent}, denoted by $p\equiv q$.

A standard high-level algorithm\footnote{
Hopcraft's celebrated algorithm \cite{MR0403320}
for minimizing a DFA has runtime complexity $O(n\log n)$.}
 for minimizing a DFA
proceeds in two stages:
\bit
\item {\tt REMOVE-UNREACHABLE}: Remove all states $q$ such that
there is no directed path from the starting state $q_0$ to $q$.
\item {\tt COLLAPSE-EQUIVALENT}: Collapse each set of mutually
equivalent states into a single state.
\eit

\section{Main results}
Our main result is the following estimate on the
state complexity of random DFAs:
\bethn
\label{thm:main}
Let $M_n\k$ be a random 
DFA 
on $n$ states and $k$ symbols
drawn uniformly from 
$[n]^{[n]\times[k]}\times2^{[n]}$.
Then, for any fixed $k\ge2$ and sufficiently large $n$,
\beqn
\label{eq:Mnkdev}
\pr{\abs{\nrm{M_n\k}-
\alpha_k n}>\sqrt n\log n
}
&=& 
\Theta\paren{\frac{
1
}{n^k}},
\eeqn
where $\alpha_k$
is unique positive root\footnote{
A closed-form expression for $\alpha_k$ 
is possible
via the Lambert $W$ function \cite{MR1809988}:
$\alpha_k = 1+W(-ke^{-k})/k$.
This constant also appears (as $\omega_k$) in \cite{DBLP:conf/stacs/BassinoDS12}
and seems to be
intimately related to
average reachability properties of semiautomata in
$[n]^{[n]\times[k]}$ under the uniform measure.
} of
$x=1-e^{-kx}$.
In particular,
\beq
\E\nrm{M_n\k} = \alpha_k n + O(\sqrt n\log n).
\eeq
\enthn
\begin{rem}
\label{rem:pittel}
Observe that 
$
0.7968\approx
\alpha_2 < \alpha_3 <\ldots<\alpha_\infty=1$.
For $k=1$, the behavior 
of $\nrm{M_n\k}$
is qualitatively different
than described in Theorem \ref{thm:main}.
The
equation $x=1-e^{-x}$
has no positive solution
and 
$\E\nrm{M_n^{(1)}} = \Theta(\sqrt n)$,
which follows from the analysis in
\cite{MR690140}.
\end{rem}
\begin{rem}
\label{rem:lb}
The lower bound of $\Omega(1/n^k)$ in (\ref{eq:Mnkdev}) is trivial,
since 
with probability 
$1/n^k$,
all of state 1's arrows point back to itself and
$\nrm{M_n\k}=1$.
\end{rem}

Our proof of Theorem \ref{thm:main} proceeds in two principal stages.
First we show that in our model, a random semiautomaton has 
roughly $\alpha_kn$ reachable states with high probability.
As in
\cite{Trakhtenbrot73},
we refer to the number of reachable states as the {\em accessibility spectrum} of the automaton.
\bethn
\label{thm:Rn}
Let
$R_n\k$ 
be
the 
accessibility spectrum 
of a 
random
semiautomaton 
on $n$ states and $k$ symbols
drawn uniformly from
$[n]^{[n]\times[k]}$.
Then, for every fixed $k\ge2$ and as $n\to\infty$,
\beqn
\label{eq:Rndev}
\pr{\abs{R_n\k-\alpha_k n}>\sqrt n\log n
}
&=& 
O\paren{\frac{1}{n^k}}.
\eeqn
\enthn

Second, we show that with high probability, very few states are lost
when equivalent ones are merged. Define $E_n\k$ to be the number
of ``excess'' reachable states: 
$$E_n\k = R_n\k - \nrm{M_n\k}.$$
Note that in principle we need only show that
the number of states lost due to merging is small {\em after} the unreachables
are removed, but we will actually show that this is true even without removing them.

\bethn
\label{thm:Enk}
For every fixed $k\ge2$,
\beq
\pr{E_n\k > 
C_k
\frac{
\log n}{\log\log n}
} = O\paren{\oo{n^k}}
\eeq
for an appropriate constant $C_k$.
\enthn
\begin{rem}
Theorem \ref{thm:Enk}
continues to hold
when each state is accepting with probability $0<p<1$ instead of $1/2$;
only the constants $C_k$ and those implicit in $O(\cdot)$ 
will change.
\end{rem}

\section{Proofs}
\belen
\label{lem:F}
Define the function
\beq
F(t) = 
n(1-(1-1/n)^t)-(t-1)/2,
\qquad 1\le t\le n
.
\eeq
Then, for sufficiently large $n$,
\beq
\frac{F^2(t)}{t} &\ge& 
\left\{
\begin{array}{lll}
0.01t,&& t\le n/2
,\\
\Omega(
\log^2n
)
,&& t\in[n/2,
\alpha_2n-\sqrt n\log n
]
\cup
[\alpha_2n+\sqrt n\log n,\infty)
.
\end{array}\right.
\eeq
\enlen
\bepf
We have $F(0)=1/2$, and for $t\le n/2$
\beq
F'(t) &=& -\oo2 +n\log\paren{1+\oo{n-1}}\cdot\paren{1-\oo n}^t \\
&\ge& -\oo2 + n\paren{\oo{n-1}-\oo{2(n-1)^2}}\cdot\paren{1-\oo n}^{n/2} \\
&\ge& -\oo2 + \frac{n}{n-1}\cdot\frac{2n-3}{2n-2}\cdot\paren{\oo{\sqrt e}-o(1)}\\
&\ge& \oo{\sqrt e}-\oo2-o(1)>0.1.
\eeq
This proves the estimate on $F^2(t)/t$ in the range $[1,n/2]$.

Now consider $t\in[n/2,
\alpha_2n-\sqrt n\log n
]$, and observe that $F(t)=H(t)+O(1)$, where
\beq
H(t)=n-t/2-n\exp(-t/n).
\eeq
By the definition of $\alpha_2$,
we have
$H(\alpha_2n)=H(0)=0$.
Furthermore, $H''(t)=-\exp(-t)/n<0$, and so $H$ is concave
with $H'(n\log2)=0$, and therefore increasing on
$[0,n\log2]$ and decreasing on $[n\log2,\infty)$. 
Hence, to lower-bound $H^2(t)/t$ in the given range, it
suffices to estimate $H$
at its right endpoint:
\beq
H(\atwo  n-\sqrt n\log n) &=& \oo2\sqrt n\log n-e^{-\atwo }\sqrt n\log n+O(\log^2 n)
=\Omega(\sqrt n\log n)
.
\eeq
Since $H'(t)<1/2-e^{-\atwo }$ for $t>\atwo n$, we have
$H(\alpha_2  n+x)=\Omega(x)$ for $x>0$,
which completes the proof.
\enpf

\bepf[Proof of Theorem \ref{thm:Rn}]
We will prove the theorem for $k=2$;  
the general case is completely analogous --- only the
constants implicit in $O(1/n^k)$ will vary with $k$.
For readability, we will write $\atwo=\alpha_2
$
and
$R_n=R_n^{(2)}$.
It will be convenient to embed 
$R_n
$ 
in a slightly more general random process.
Fix $n\ge1$, and
define the 
sequence of
random variables $(\nu_t)_{t=1}^\infty$,
as follows:
\beq
\nu_1 &=& 1,\\
\nu_{t+1} &=&
\left\{
\begin{array}{ll}
\nu_t,& \text{with probability } \nu_t/n
,\\
\nu_t+1,& \text{with probability } 1-\nu_t/n.
\end{array}\right.
\eeq
Clearly, 
$\nu_t$ is with probability 1 nondecreasing,
upper-bounded by $n$, and reaches $n$ after a finite
number of steps.
Let us also define
\beqn
\label{eq:ON}
\omega_t
=
2\nu_t
+ 
1-t,
\qquad t\ge1.
\eeqn

Now consider the
following process for generating random 
directed multigraphs with regular outdegree 2.
For time steps $t=1,2,\ldots$,
we will maintain the set of nodes $N_t$,
reached from $q_0=1$ by time $t$,
and two sets of edges: {\em open edges} $O_t$
and {\em closed edges} $C_t$.
A closed edge $c$ is an ordinary directed arrow from a source node $p$ to a 
destination node $q$ marked with a $\sigma\in[k]$ and
denoted by $c=(p\wed{\sigma}q)$.
An open edge $o$ has a specified source $p$ but an as yet unspecified destination;
such an edge will be denoted by $o=(p\wed{\sigma}\star)$.
We initialize $N_1=\set{1}$,
$C_1=\emptyset$
and $O_1=\set{(1\wed{0}\star),(1\wed{1}\star)}$.
At time $t+1$, some (arbitrarily chosen\footnote{It is easy to see
that the distribution of $N_t,C_t,O_t$ is unaffected by the order
in which the open edges are selected.}) 
open edge in $o
\in O_t$ 
(if one exists)
chooses a destination node $q$
as follows:
\bit
\item[(i)] 
$q\in N_t$
with probability $|N_t|/n$
(that is, $o$ will point to a previously reached node);
\item[(ii)]
$q\in[n]\setminus N_t$
with probability $1-|N_t|/n$.
\eit
In event (i), $O_{t+1}=O_t\setminus\set{o}$,
while
in event (ii), $O_{t+1}=(O_t\setminus\set{o})\cup
\set{(q\wed{0}\star),(q\wed{1}\star)}$;
in both cases, $N_{t+1}=N_t\cup\set{q}$
and $C_{t+1}=C_t\cup\set{o}$.

The random semiautomaton embeds into the process $(\nu_t,\omega_t)$  
via the following natural
correspondence: $|N_t|=\nu_t$ and $|O_t|=\omega_t$ as long as
 the latter is nonnegative (in particular, 
the correspondence breaks down 
for $t>2n+1$,
since 
$\omega_t$ becomes negative).
Let $\tau$ be
the smallest $t$ for which $\omega_t=0$
--- 
i.e., the first time
there are no longer any open edges to choose from.
Then the pair $(N_\tau,C_\tau)$ defines\footnote{
Since the quantity of interest is the 
accessibility spectrum,
it is unnecessary to define transitions out of unreachable states.
}
a semiautomaton with 
accessibility spectrum
$R_n=\nu_\tau$, drawn uniformly from 
$[n]^{[n]\times\set{0,1}}$.
Hence, proving (\ref{eq:Rndev}) amounts to showing that
\beqn
\label{eq:taudev}
\pr{\abs{\tau-2\atwo  n}>
\sqrt n\log n
} &=& O\paren{\frac{
1
}{n^2}}.
\eeqn
Indeed, (\ref{eq:taudev}) implies that $\tau=(2+o_P(1))\atwo  n$,
and $\nu_\tau=(1+o_P(1))\nu_{\atwo n}$.
Since,
by definition, $\tau$ is the smallest $t$ for which
$\nu_t=(t-1)/2$,
we have
\beqn
\nonumber
\pr{\tau\in[a,b]}
&\le&
\pr{\exists t\in[a,b
]: \nu_t=(t-1)/2} \\
\label{eq:tau[a,b]}
&\le&
\pr{\exists t\in[a,b
]: \nu_t\le(t-1)/2}
.
\eeqn
We estimate the left tail of $\tau$ as follows:
\beq
\pr{\tau\le\atwo n-
\sqrt n\log n
} &\le&
P_0 + P_1 + P_2,
\eeq
where
\beq
P_0 &=& \pr{\tau\in[1,150\log n]},\\
P_1 &=& \pr{\tau\in[150\log n,n/2]},\\
P_2 &=& \pr{\tau\in[n/2,\atwo n-
\sqrt n\log n
]}.
\eeq
To bound $P_0$, we
first argue, by elementary combinatorics, that $\pr{\omega_3<3}=O(1/n^2)$.
Now we condition on the high-probability event that there are at least $3$ open arrows available
after $3$ steps. 
If all of the open arrows have been exhausted between time $t=4$ and $t=T$, then certainly at least
three of these arrows must point back to the $O(T)$ previous states. Thus, for $T=150\log n$,
\beq
P_0 \in 
O\paren{\frac{1}{n^2}
+
{\binom{T}{3}}\paren{\frac{T}{n}}^3} 
\subset 
O\paren{\frac{1}{n^2}}
.
\eeq

To bound $P_1$ and $P_2$, we observe that an alternate interpretation
is possible for $\nu_t$.
Namely, 
when $t$ balls are thrown into $n$ bins 
uniformly at random,
the number of non-empty bins 
is distributed as $\nu_t$.
We also observe that
\beqn
\label{eq:Enut}
\E\nu_t = n(1-(1-1/n)^{t})
,
\qquad t\ge1
.
\eeqn
Now by the Chernoff bound for negatively associated random variables
\cite[Prop. 5, Thm. 13]{Dubhashi:1998:BBS:299633.299634}, 
\beqn
\label{eq:duvshani}
\pr{\nu_t-\E\nu_t \le 
-\Delta
}
&\le& \exp(-2\Delta^2/t),
\qquad \Delta>0.
\eeqn
Hence,
\beq
P_1 &\le& \sum_{t=150\log n}^{n/2} \pr{\nu_t-\E\nu_t\le -F(t)
},\\
P_2 &\le& \sum_{t=n/2}^{\atwo n-
\sqrt n\log n
} \pr{\nu_t-\E\nu_t\le -F(t)},
\eeq
where $F(t)$
is as in Lemma \ref{lem:F}.
The estimates in the lemma and (\ref{eq:duvshani})
yield
\beq
P_1 &\le& \frac{n}{2}\exp(-3\log n) \in O\paren{\oo{n^2}}
\eeq
and
\beq
P_2 &\in& O(n)\exp(-\Omega(\log^2 n)) 
\subset O\paren{\oo{n^2}}.
\eeq
This proves the left-tail estimate in (\ref{eq:taudev}).
To prove the corresponding right-tail estimate,
we observe that, analogously to (\ref{eq:tau[a,b]}),
\beq
\pr{\tau\in[a,b]}
\le
\pr{\exists t\in[a,b
]: \nu_t\ge(t-1)/2}
.
\eeq
The deviation probability is bounded as in (\ref{eq:duvshani}):
\beq
\pr{\nu_t-\E\nu_t\ge
\Delta
}
&\le& \exp(-2\Delta^2/t),
\qquad \Delta>0.
\eeq
Hence
\beq
\pr{\tau>\atwo n+\sqrt n\log n
} &\le& \sum_{t=\atwo n+\sqrt n\log n}^n \pr{\nu_t-\E\nu_t\ge G(t)},
\eeq
where $G(t)=-F(t)$.
Invoking again Lemma \ref{lem:F}, we have
\beq
\pr{\tau>\atwo n+\sqrt n\log n
} &\in& 
O(n)
\exp(-\Omega(\log^2 n)) 
\subset O\paren{\oo{n^2}}.
\eeq
\enpf

\bepf[Proof of Theorem \ref{thm:Enk}]
Again, for
ease of exposition, we only prove the claim for $k=2$.
We start by explaining the idea of the proof. We need to show that there are usually ``few'' pairs of equivalent states.
Let us start by describing two ``typical'' situations in which equivalent states emerge.
The first is where a state is mapped into itself by every member of $\set{0,1}$. 
Two such states are equivalent if and only if both are accepting or both are rejecting, which happens with a probability
of $1/2$. More generally, if from each of the two states one can reach very few states, then there is a non-negligible probability
that the states are equivalent. Thus, we will show that there are few states 
with small accessibility spectra.
In the preceding sentence,
``few'' means (with high probability) ``at most $2$'', while ``small'' means ``less than $4\log_2n$''.

The second principal reason for two states $q,q'$ to be equivalent is that 
$\delta(q,0)=\delta(q',0)$
and
$\delta(q,1)=\delta(q',1)$.
Again, $q$ and $q'$ are equivalent in this case with probability $1/2$.
Thus, we will need to show that there are few pairs of states $q,q'$
for which there are few words in $\set{0,1}^*$ taking $q$ and $q'$ to distinct states.
Here, the first ``few'' means ``at most $C\log n/\log\log n$'' and the second means ``up to $4\log_2 n$''.

Let us now consider the above scenarios in more detail. 
The (random) set of states reachable from $q$ is given by $\set{q,\delta(q,0),\delta(q,1),\delta(q,00),\delta(q,01),\ldots}$.
Thus, the states reachable from $q$ reside on a binary tree whose edges are marked by letters in $\set{0,1}$.
Each time the random DFA selects a state $p=\delta(q,w)$, if $p$ is already in the tree, the edge that would create a directed cycle is not drawn.
We refer to the resulting tree as the tree {\em growing} from $q$.
Its size is the {accessibility spectrum} of $q$, denoted by $S(q)$.

Let $C>0$ be a constant to be determined later.
A state's accessibility spectrum is said to be {\em small} if it is below $C\log_2n$.
As in the
proof of Theorem \ref{thm:Rn},
the probability of a given state having a small accessibility spectrum is $O(1/n^2)$.
A similar argument shows that the joint probability of any pair of states
$q,q'$ 
having small accessibility spectra is $\tilde O(1/n^4)$. 
Indeed, 
consider the event of $S(q')$ being small, conditioned on $S(q)$ being such.
Draw the states ${\delta(q',0),\delta(q',1),\delta(q',00),\ldots}$ similarly 
to the proof of Theorem \ref{thm:Rn}. The event in question is contained in
the event whereby, in the course of the first $C\log_2n$ steps of the process
of ``closing'' the open edges, we encounter at least twice
either a state visited already or a state belonging to the tree growing from $q$.
The probability of the latter event is clearly $O(\log^4/n^2)$.
Hence,
\beq
\pr{S(q),S(q')\text{ are both small}}
=O\paren{\frac{\log^4n}{n^4}}.
\eeq
Carrying this line of reasoning over to triples,
we have that the probability of any
three states having small accessibility spectra is $\tilde O(1/n^6)$
--- and therefore,
\beq
\pr{\text{there are 3 distinct states with small accessibility spectra}}
&\in &
\tilde O\paren{\oo{n^6}
\binom{n}{3}}\\
&=&\tilde O\paren{\oo{n^3}}
\subset O\paren{\oo{n^2}}.
\eeq

In view of the discussion above, we may assume (after removing up to 2 states)
that all states have large accessibility spectra. Consider two states $q,q'$. Let $T$ be a tree
of size $m=C\log_2n$ growing from $q$ (this will typically be a subtree of a larger tree of size $O(\alpha n)$). 
The nodes of $T$ are ${\delta(q,w_1),\delta(q,w_2),\ldots,\delta(q,w_m)}$ for certain words 
$w_1,w_2,\ldots,w_m\in\set{0,1}^*$.
If $\delta(q,w_i)\neq\delta(q',w_i)$,
$i=1,2,\ldots,m$, then the probability of $q,q'$ being equivalent is at most $1/2^m=1/n^C$.
(Note that this holds even if the states $\delta(q',w_i)$, $1\le i\le m$, are not
mutually distinct, in fact even if they all coincide. Similarly,
it does not matter if some of the states $\delta(q',w_i)$ coincide with some of the
$\delta(q,w_j)$, as long as $i\neq j$.)
The probability that both $\delta(q,0)=\delta(q',0)$
and
$\delta(q,1)=\delta(q',1)$ is $1/n^2$. 
Call a state pair satisfying these equalities a {\em dud}.
The union bound does not yield a non-trivial upper bound on the number of duds,
and a more refined analysis will be needed. 
Clearly, the probability that $d$ 
specific 
pairwise disjoint
pairs are duds is $1/n^{2d}$.
Now the probability that there exist $d$ disjoint duds is at most
\beq
\label{eq:duds}
\oo{n^{2d}}\binom{n}{2d}\cdot(2d-1)\cdot(2d-3)\cdot\ldots\cdot1
\le
\oo{n^{2d}}\cdot \frac{n^{2d}}{(2d)!}\cdot\frac{(2d)!}{2^dd!}=\oo{2^dd!}
.
\eeq
Choosing $d=3\log n/\log\log n$ 
and applying Stirling's formula, we see that
the probability 
of there being $d$ disjoint duds is $O(1/n^2)$.

Other than duds --- pairs ``dying'' right away after 2 steps ---
we must consider pairs dying after $4, 6,\ldots,C\log_2n$ steps. However,
the probability that a pair will die after 4 steps is $O(1/n^3)$, after
6 steps --- $O(1/n^4)$, and so forth. Hence, the probability
that there will be 2 pairs for which the process dies after 4 steps is
\beq
\binom{n}{4}\cdot O\paren{\frac{1}{n^6}}= O\paren{\oo{n^2}},
\eeq
that there will be a pair that dies after 6 steps,
\beq
\binom{n}{2}\cdot O\paren{\frac{1}{n^4}}= O\paren{\oo{n^2}},
\eeq
and that some pair will die after 
$t\in[8,4\log_2n]$ steps,
\beq
\binom{n}{2}\cdot\tilde O\paren{\frac{1}{n^5}}= \tilde O\paren{\oo{n^3}}\subset O\paren{\oo{n^2}}.
\eeq
Now for two states $q,q'$
reaching distinct states for 
many words $w_i$,
the probability of being equivalent is at most $1/n^C$.
Thus, it suffices to take $C=4$
to bound the probability of any pair of states
growing large trees yet being equivalent by
\beq
\binom{n}{2}\cdot O\paren{\frac{1}{n^4}}= O\paren{\oo{n^2}}.
\eeq
\enpf

\bepf[Proof of Theorem \ref{thm:main}]
Follows immediately from Theorems \ref{thm:Rn} and \ref{thm:Enk}
since the former estimates the number of 
states remaining after
{\tt REMOVE-UNREACHABLE}
and the latter bounds the number of states lost after
{\tt COLLAPSE-EQUIVALENT}.
\enpf

\bibliographystyle{plain}
\bibliography{../mybib}

\end{document}